\documentclass{notices}
\usepackage{amsfonts,amssymb,amsmath,amscd}
\usepackage{mathtools}
\usepackage{tikz}

\title{
What is a Parking Function?
}

\author{
    J. Carlos Mart\'inez Mori\affil{
        J. Carlos Mart\'inez Mori is a Schmidt Science Fellow and a President's Postdoctoral Fellow at the Georgia Institute of Technology.
        His email address is jcmm@gatech.edu.
    }
}

\begin{document}

\maketitle

\emph{Abstract}.
In this expository article I describe classical results in the combinatorics of parking functions.
Its English-Spanish translation is included.
\\

Consider a one-way street with $n \in \{1, 2, \ldots\}$ numbered parking spots, denoted $s_1, s_2, \ldots, s_n$.
A sequence of $n$ cars enter the street one at a time, each with a preferred spot.
Upon its arrival, car $1 \leq i \leq n$, denoted $c_i$, drives to its preferred spot $x_i \in [n] \coloneqq \{1, 2, \ldots, n\}$.
If spot $x_i$ is unoccupied, car $c_i$ is \emph{lucky} and parks there.
Otherwise, car $c_i$ \emph{displaces} further down the street until it finds the first unoccupied spot in which to park, if such a spot exists.
If no such spot exists, car $c_i$ reneges the search process unable to park.
Let $\mathbf{x} = (x_1, x_2, \ldots, x_n) \in [n]^n$ be the $n$-tuple encoding the cars' parking preferences.
If all cars are able to park, then $\mathbf{x}$ is said to be a \emph{parking function} of length $n$. 

\begin{figure}[ht]
    \centering
    \includegraphics[clip, trim=31cm 19cm 31cm 19cm, width=\linewidth]{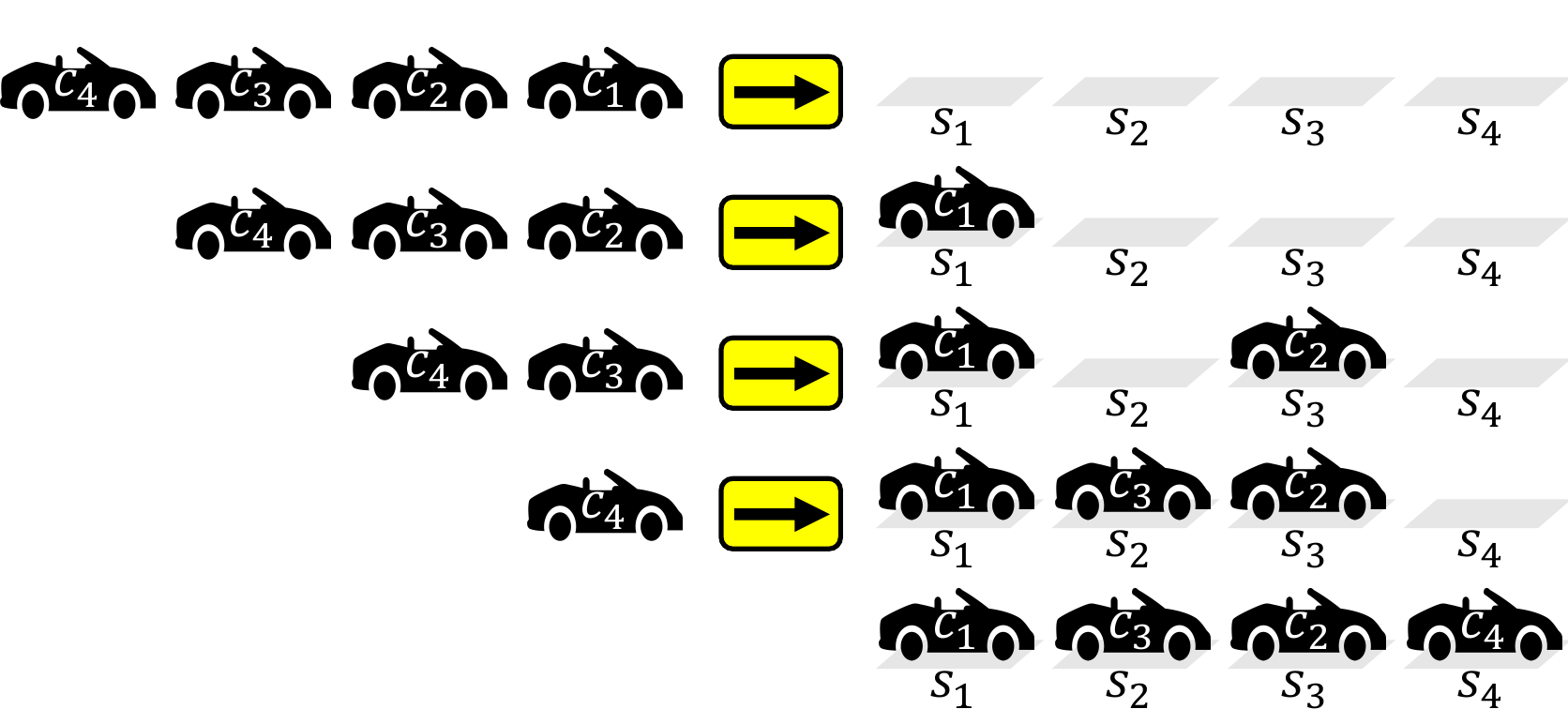}
    \caption{
        The parking function $(1,3,1,1)$ has outcome $(1,3,2,4)$.   
        Cars $c_1$ and $c_2$ are lucky and park in their preferred spots $s_1$ and $s_3$, whereas cars $c_3$ and $c_4$ displace further down the street before ultimately parking in spots $s_2$ and $s_4$, respectively.
    }
    \label{fig: parking example}
\end{figure}

Parking functions were first implicitly studied by Pyke~\cite{pyke1959supremum} in his study of Poisson processes and later on by Konheim and Weiss~\cite{konheim1966occupancy} in their study of hashing with linear probing.
Parking functions can be seen as functions in the sense that the ``parking experiment'' narrated above assigns, to each $n$-tuple of parking preferences, a single $n$-tuple encoding its parking \emph{outcome}.
If the $n$-tuple of preferences is indeed a parking function, its outcome can be treated as a permutation of $[n]$ written in one-line notation.
For example, as depicted in Figure~\ref{fig: parking example}, $(1,3,1,1)$ is a parking function of length $4$ with outcome $(1,3,2,4)$.

Classical enumerative results about parking functions foreshadow their rich mathematical structure.
There are
\begin{equation}
\label{eq: cayley}
    (n+1)^{n-1}
\end{equation}
parking functions of length $n \geq 1$ (OEIS \href{https://oeis.org/A000272}{A000272})\textemdash this is Cayley's formula~\cite{cayley1889theorem} for the number of labeled trees on $n + 1$ vertices.
This count was established independently by Pyke~\cite{pyke1959supremum} and Konheim and Weiss~\cite{konheim1966occupancy}, and has been recovered bijectively by numerous authors; refer to Yan~\cite{yan2015parking}*{Section 13.2} and the references therein.

\begin{figure}[ht]
    \centering
    \includegraphics[clip, trim=32cm 15cm 32cm 10cm, width=\linewidth]{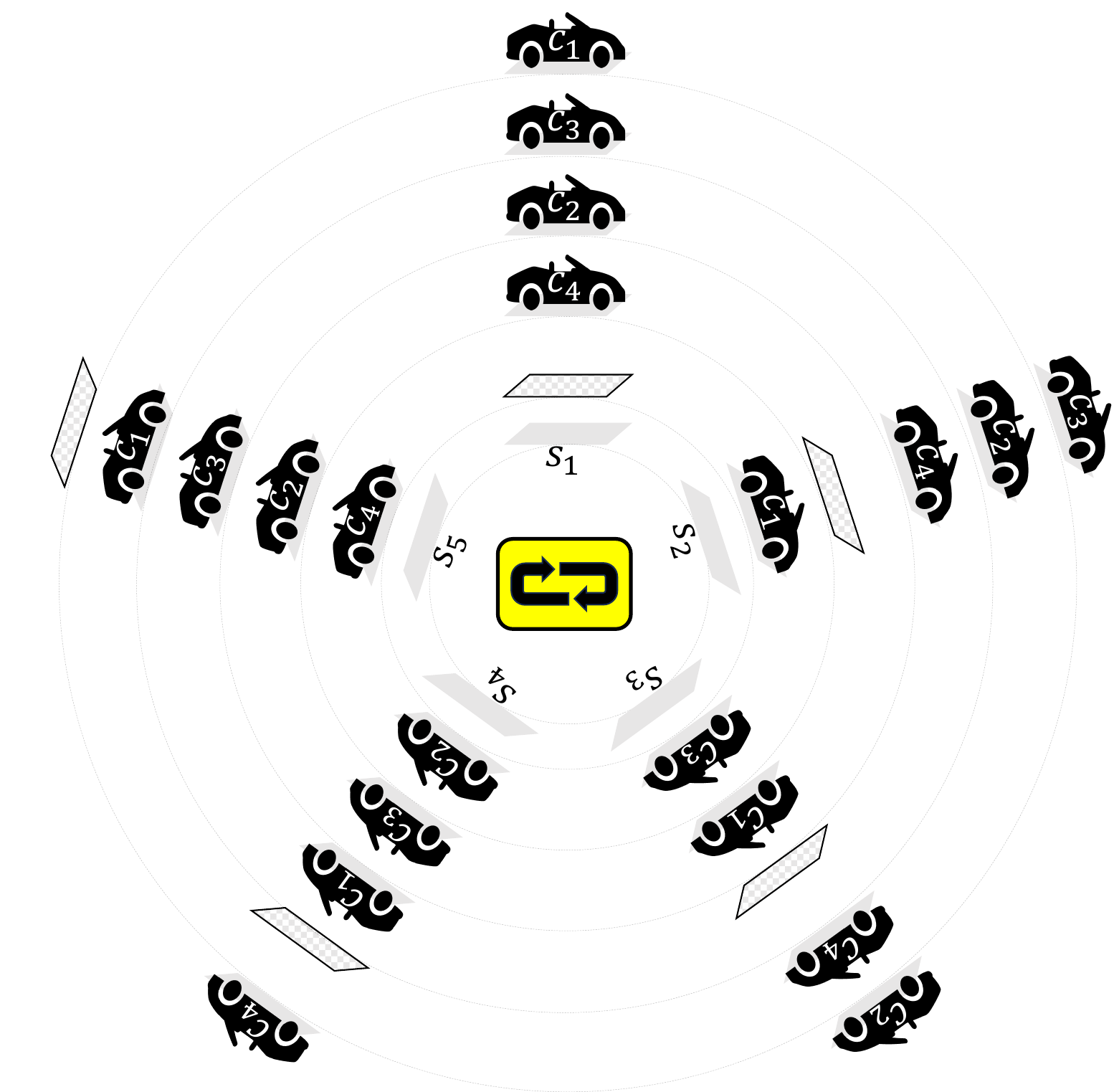}    
    \caption{
        Pollak's argument for \eqref{eq: cayley} using the example in Figure~\ref{fig: parking example}.
        The tuples $(1,3,1,1)$, $(5,2,5,5)$, $(4,1,4,4)$, $(3,5,3,3)$, and $(2,4,2,2)$ form an equivalence class of order $5$.
        Their outcomes are depicted as layers overlaying a circular one-way road with spots $s_1, s_2, \ldots, s_5$.
        The only parking function in this class is $(1,3,1,1)$ because it is the only tuple that leaves spot $s_5$ unoccupied in its outcome, as depicted in the outermost layer of the illustration.
    }
    \label{fig: pollak}
\end{figure}

An elegant proof of \eqref{eq: cayley}, due to Pollak as credited by Riordan~\cite{riordan1969ballots}, is as follows.
Consider a circular one-way street with $n + 1$ numbered parking spots, denoted $s_1, s_2, \ldots, s_n, s_{n+1}$.
Say, it is a roundabout with no exits and a single entry between spots $s_{n+1}$ and $s_1$.
There are $n$ cars attempting to park, so that each $n$-tuple $\mathbf{x} = (x_1, x_2, \ldots, x_n) \in [n+1]^n$ encodes a possibility for the cars' parking preferences.
How many of these $n$-tuples are in fact parking functions of length $n$?
Since there are $n+1$ spots but only $n$ cars, there will always be one spot left unoccupied; which exact spot this is depends on the preferences.
Note that the preference tuples form a group $G$ under component-wise addition modulo $n + 1$, and consider the subgroup $H \leq G$ generated by the all-ones preference tuple $(1, 1, \ldots, 1)$.
As depicted in Figure~\ref{fig: pollak}, the cosets $\mathbf{x}H = \{\mathbf{x} + h : h \in H\}$ for $\mathbf{x} \in G$ form equivalence classes of order $n + 1$, each of which contains exactly one parking function, namely whichever tuple leaves spot $n + 1$ unoccupied in its outcome.
Therefore, there are $(n+1)^n / (n+1) = (n+1)^{n-1}$ parking functions of length $n$.

Now, consider an arbitrary $n$-tuple $\mathbf{x} = (x_1, x_2, \ldots, x_n) \in [n]^n$.
Can one tell whether it is a parking function of length $n$ without explicitly considering its parking outcome?
A classical inequality-based characterization of parking functions allows precisely this.
Let $\mathbf{x}' = (x_1', x_2', \ldots, x_n')$ be the weakly increasing rearrangement of $\mathbf{x}$.
Then, the following holds:
\begin{align}
\label{eq: characterization}
    \parbox{0.375\linewidth}{
        \centering 
        $\mathbf{x}$ is a parking function of length $n$
    } 
    &&
    \iff
    &
    \parbox{0.375\linewidth}{
        \centering 
        $x_i' \leq i$ for all $1 \leq i \leq n$
    }
\end{align}
This characterization can be verified as follows.
If there exists some $i$ for which $x_i' > i$, then at least $n - i + 1$ cars attempt to park in the last $n - i$ spots $s_{i + 1}, s_{i + 2}, \ldots, s_n$.
This is more cars than spots available, so at least one of these cars will be unable to park.
Conversely, if some car is unable to park, then there exists an earliest spot $s_i$ left unoccupied, which can only be the case if $x_i' > i$.

Note that \eqref{eq: characterization} has the following remarkable implication: parking functions are invariant under the action of the symmetric group $\mathfrak{S}_n$, which permutes their subscripts.
In particular, the set of parking functions is obtained from the rearrangements of weakly increasing parking functions.
For example, $(1,3,1,1)$, $(1,1,3,1)$, $(1,1,1,3)$, and $(3,1,1,1)$ are all parking functions of length $4$, each with a different parking outcome, whereas no rearrangement of $(1,3,3,4)$ can possibly be a parking function of length $4$.

Parking functions are also closely related to the Catalan numbers (OEIS \href{https://oeis.org/A000108}{A000108}), which are given by the recurrence relation
\begin{equation}
\label{eq: catalan}
    C_{n} 
    = \sum_{i=1}^{n} C_{i-1} \cdot C_{n - i}
\end{equation}
for $n \geq 1$ with $C_0 = 1$.
In particular, the set of weakly increasing parking functions of length $n$ is enumerated by \eqref{eq: catalan}.
To verify this, suppose you construct a weakly increasing parking function of length $n$, denoted $\mathbf{x}' = (x_1', x_2', \ldots, x_n') \in [n]^n$, with your choice of $1 \leq i \leq n$ for the greatest index satisfying $x_i' = i$.
For any such choice of $i$, the $(i-1)$-tuple $(x_1', x_2', \ldots, x_{i-1}')$ must be a parking function of length $i-1$, the $(n-i)$-tuple $(x_{i+1}' - i + 1, x_{i+2}' - i + 1, \ldots, x_n' - i + 1)$ must be a parking function of length $n-i$, and by induction there are $C_{i-1} \cdot C_{n - i}$ possibilities.

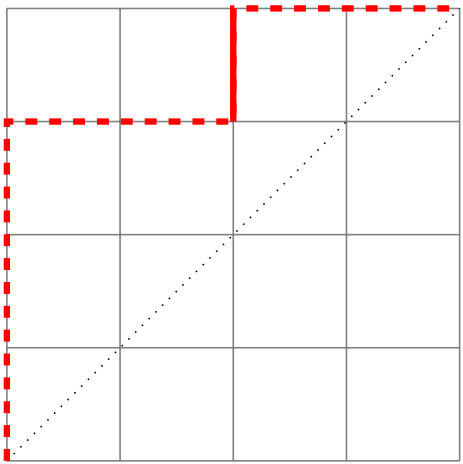
\begin{figure}[ht]
\centering
\resizebox{0.8\linewidth}{!}{
\begin{tikzpicture}
    \draw[step=1cm,gray,thin] (0,0) grid (4,4);
    \draw[thin, dotted] (0,0)--(4,4);
    \draw[ultra thick,dashed,red] (0,0)--(0,1)--(0,2)--(0,3)--(1,3)--(2,3)--(2,4)--(3,4)--(4,4);
    \draw[ultra thick,red] (2,3)--(2,4);
\end{tikzpicture}
}
    \caption{
        The Dyck path $NNNEENEE$ of length $8$, illustrated with bold dashed red lines, maps to the weakly increasing parking function $\mathbf{x}' = (x_1', x_2', x_3', x_4') = (1,1,1,3)$ of length $4$.
        For example, the fourth $N$ step in the path is the sixth step overall, as illustrated with a solid bold red line.
        There are a total of two $E$ steps before it, so $x_4' = 1 + 2 = 3$.
    }
    \label{fig: dyck path}
\end{figure}

As a consequence, weakly increasing parking functions are in bijection with the wide variety of Catalan objects; refer to Stanley~\cite{stanley2015catalan} for a survey.
For example, weakly increasing parking functions of length $n$ are in bijection with Dyck paths of length $2n$; these are lattice paths from $(0,0)$ to $(n,n)$ using only north steps $(1,0)$ and east steps $(0,1)$, denoted ``N'' and ``E'' respectively, and which do not cross below the main diagonal.
Armstrong, Loehr, and Warrington~\cite{armstrong2016rational}*{Section~2.2} describe the following bijection: given a Dyck path of length $2n$, obtain a weakly increasing parking function $\mathbf{x}' = (x_1', x_2', \ldots, x_n')$ of length $n$ by letting $x_i'$ be one plus the total number of E steps appearing before the $i$th N step for all $1 \leq i \leq n$.
Figure~\ref{fig: dyck path} illustrates this construction using the weakly increasing parking function $(1,1,1,3)$; the weakly increasing rearrangement of the example in Figure~\ref{fig: parking example}.

Note that if $\mathbf{x}' = (x_1', x_2', \ldots, x_n') \in [n]^n$ is a weakly increasing parking function, then the $i$th car $c_i$ (with preference $x_i'$) parks in the $i$th spot $s_i$.
Therefore, the total displacement of a weakly increasing parking function $\mathbf{x}'$ is given by
\begin{equation}
\label{eq: displacement}
    \sum_{i=1}^n i - x_i' = \frac{n(n+1)}{2} - \sum_{i=1}^n x_i'.
\end{equation}
This statistic has further combinatorial interpretations.
For example, the number of full squares between a Dyck path and the main diagonal, which in the case of Figure~\ref{fig: dyck path} is $4$, is the same as the total displacement of its corresponding weakly increasing parking function.
Lastly, note that the dependency of \eqref{eq: displacement} on $\mathbf{x}'$ reduces to the summation of its terms.
This implies that the total displacement of a parking function is also preserved under rearrangements!

\begin{figure}[ht]
    \centering
    \includegraphics[width=\linewidth]{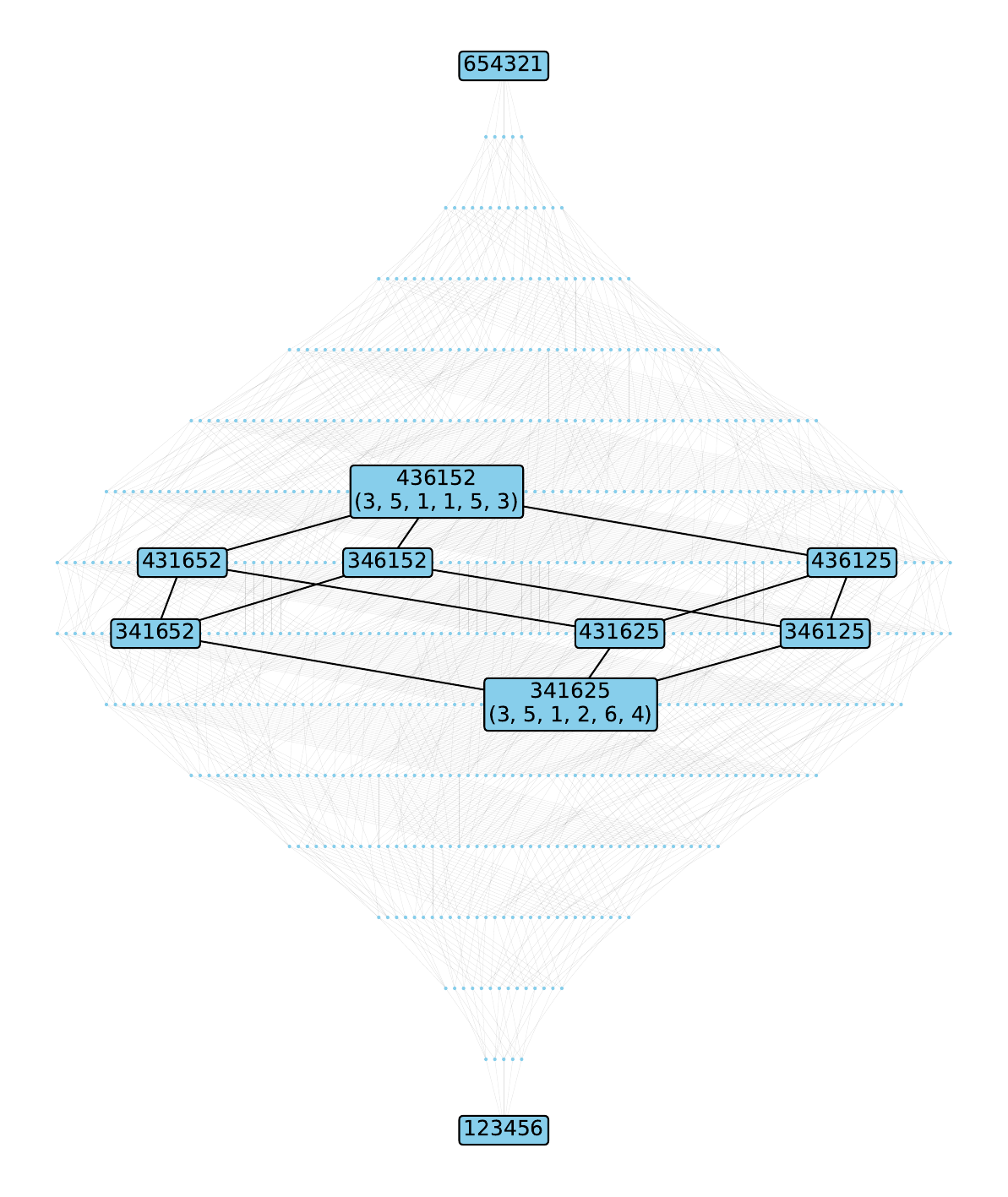}   
    \caption{
        Weak order poset of $\mathfrak{S}_6$ with a Boolean interval of rank $3$ highlighted, with minimal element $341625$ and maximal element $436152$ written in one-line notation.
        The outcome of $(3,5,1,2,6,4)$, a parking function of length $6$, is $(3, 4, 1, 6, 2, 5)$ and corresponds to the minimal element of the interval.
        The outcome of $(3,5,1,1,5,3)$, another parking function
        of length $6$, is again $(3, 4, 1, 6, 2, 5)$ and corresponds to the minimal element of the interval.
        In fact, $(3,5,1,1,5,3)$ is a ``unit Fubini ranking with $3$ distinct ranks'' and corresponds to a Boolean interval of rank $6-3 = 3$ (i.e., a cube), whereas $(3,5,1,2,6,4)$ is a ``unit Fubini ranking with $6$ distinct ranks'' and corresponds to a Boolean interval of rank $6 - 6=0$ (i.e., a node).
        Refer to Elder et al.~\cite{elder2024boolean} for details.
    }
    \label{fig: weak}
\end{figure}

Despite first appearing in the literature more than six decades ago, the combinatorics of parking functions remains a vibrant area of research.
Recent work ranges from the study of discrete statistics such as displacement~\cites{kenyon2023parking,elder2023cost} or the number of lucky cars~\cites{gessel2006refinement,slavik2023lucky,stanley2023some}, the introduction of new variants and/or generalizations on the classical ``parking experiment'' rule (refer to Carlson et al.~\cite{carlson2021parking} for an accessible tour of endless possibilities in the style of ``choose your own adventure''), polyhedral aspects~\cites{amanbayeva2022convex,hanada2023generalized}, and surprising connections to seemingly unrelated objects~\cites{aguillon2023parking,harris2024lucky}, to mention just a few.
In recent work, my collaborators and I use a subset of parking functions we call \emph{unit Fubini rankings}, and in particular their outcome map, to characterize and enumerate the Boolean intervals of rank $k$ in the weak order poset of $\mathfrak{S}_n$~\cite{elder2024boolean}.
Figure~\ref{fig: weak} illustrates our construction.

\bibliography{bib}

\end{document}


\maketitle

\emph{Resumen}.
En este artículo de difusión matemática describo resultados clásicos en la combinatoria de funciones de parqueo.
Su traducción español-inglés está incluida.
\\

Considere una calle de una sola vía con $n \in \{1, 2, \ldots\}$ puestos de parqueo numerados, denotados $s_1, s_2, \ldots, s_n$.
Una secuencia de $n$ carros ingresa a la calle uno a la vez, cada cual con su puesto preferido.
Tras la llegada del carro $1 \leq i \leq n$, denotado $c_i$, este se dirige a su puesto preferido $x_i \in [n] \coloneqq \{1, 2, \ldots, n\}$.
Si el puesto $x_i$ está libre, el carro $c_i$ es \emph{suertudo} y se parquea ahí.
De lo contrario, el carro $c_i$ se \emph{desplaza} más a lo largo de la calle hasta encontrar el primer puesto libre en el que parquear, si tal puesto existe.
Si tal puesto no existe, el carro $c_i$ abandona el proceso de búsqueda sin poder parquear.
Sea $\mathbf{x} = (x_1, x_2, \ldots, x_n) \in [n]^n$ la $n$-tupla que codifica las preferencias de parqueo de los carros.
Si todos los carros logran parquear, se dice que $\mathbf{x}$ es una \emph{función de parqueo} de longitud $n$. 

\begin{figure}[ht]
    \centering
    \includegraphics[clip, trim=31cm 19cm 31cm 19cm, width=\linewidth]{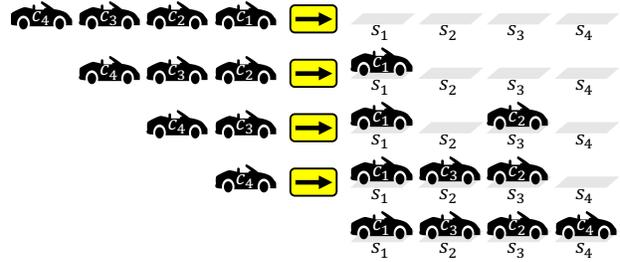}
    \caption{
        La función de parqueo
        $(1,3,1,1)$ tiene desenlace $(1,3,2,4)$. 
        Los carros $c_1$ y $c_2$ son suertudos y se parquean en sus puestos preferidos $s_1$ y $s_3$, mientras que los carros $c_3$ y $c_4$ se desplazan más a lo largo de la calle hasta finalmente parquear en los puestos $s_2$ y $s_4$, respectivamente.
    }
    \label{fig: parking example}
\end{figure}

Las funciones de parqueo fueron inicialmente estudiadas de manera implícita por Pyke~\cite{pyke1959supremum} en su estudio sobre procesos de Poisson y luego por Konheim y Weiss~\cite{konheim1966occupancy} en su estudio sobre funciones resumen con sondeo lineal.
Las funciones de parqueo se pueden interpretar como funciones de manera que el ``experimento de parqueo'' relatado anteriormente asigna, a cada $n$-tupla de preferencias de parqueo, una $n$-tupla única que codifica su \emph{desenlace} de parqueo.
Si la $n$-tupla de preferencias de parqueo es efectivamente una función de parqueo, su desenlace se puede considerar como una permutación de $[n]$ escrita en notación de una línea.
Por ejemplo, tal como descrito en la Figura~\ref{fig: parking example}, $(1,3,1,1)$ es una función de parqueo de longitud $4$ con desenlace $(1,3,2,4)$.

Los resultados enumerativos clásicos acerca de las funciones de parqueo presagian su abundante estructura matemática.
Hay
\begin{equation}
\label{eq: cayley}
    (n+1)^{n-1}
\end{equation}
funciones de parqueo de longitud $n \geq 1$ (OEIS \href{https://oeis.org/A000272}{A000272})\textemdash esta es la fórmula de Cayley~\cite{cayley1889theorem} para el número de árboles etiquetados con $n + 1$ vértices.
Este conteo fue establecido de manera independiente por Pyke~\cite{pyke1959supremum} al igual que por Konheim y Weiss~\cite{konheim1966occupancy}, y ha sido restablecido de manera biyectiva por varios autores; consulte Yan~\cite{yan2015parking}*{Sección 13.2} y las referencias ahí encontradas.

\begin{figure}[ht]
    \centering
    \includegraphics[clip, trim=32cm 15cm 32cm 10cm, width=\linewidth]{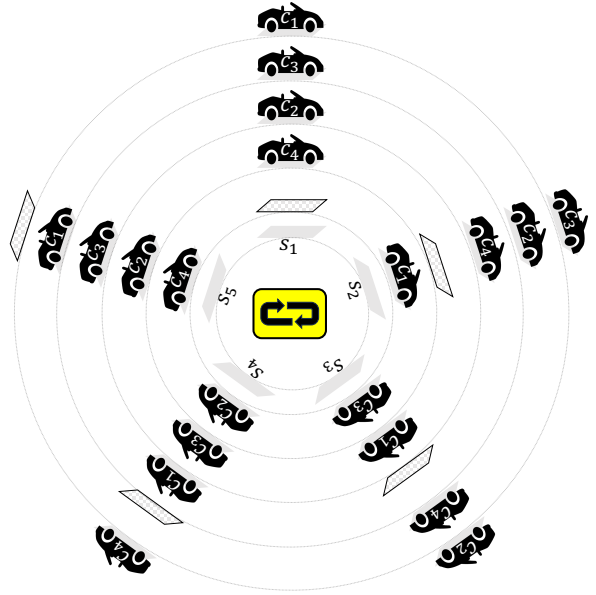}    
    \caption{
        El argumento de Pollak para \eqref{eq: cayley} utilizando el ejemplo de la Figura~\ref{fig: parking example}.
        Las tuplas $(1,3,1,1)$, $(5,2,5,5)$, $(4,1,4,4)$, $(3,5,3,3)$, y $(2,4,2,2)$ forman una clase de equivalencia de orden $5$.
        Sus desenlaces son representados mediante capas superimpuestas sobre la calle circular de una sola vía con puestos $s_1, s_2, \ldots, s_5$.
        La única función de parqueo en esta clase es $(1,3,1,1)$ ya que es la única tupla que deja el puesto $s_5$ desocupado en su desenlace, tal como reflejado en la capa exterior de la ilustración.
    }
    \label{fig: pollak}
\end{figure}

Una elegande demostración de \eqref{eq: cayley}, de parte de Pollak según Riordan~\cite{riordan1969ballots}, es la siguiente.
Considere una calle circular de una sola vía con $n + 1$ puestos de parqueo numerados, denotados $s_1, s_2, \ldots, s_n, s_{n+1}$.
Digamos que es un redondel sin salida y con una sola entrada entre los puestos $s_{n+1}$ y $s_1$.
Hay $n$ carros intentando parquear, de tal forma que cada $n$-tupla $\mathbf{x} = (x_1, x_2, \ldots, x_n) \in [n+1]^n$ codifica una posibilidad con respecto a las preferencias de parqueo de los carros.
?`Cuántas de estas $n$-tuplas son efectivamente funciones de parqueo de longitud $n$?
Como hay $n+1$ puestos pero solamente $n$ carros, siempre quedará un puesto sin ocupar; qué puesto este es exactamente depende de las preferencias.
Note que las tuplas de preferencias forman un grupo $G$ bajo suma por componentes modulo $n + 1$, y considere el subgrupo $H \leq G$ generado por la tupla de preferencias todas uno $(1, 1, \ldots, 1)$.
Tal como descrito en la Figura~\ref{fig: pollak}, las clase laterales $\mathbf{x}H = \{\mathbf{x} + h : h \in H\}$ donde $\mathbf{x} \in G$ forman clases de equivalencia de orden $n + 1$, cada cual conteniendo exactamente una función de parqueo, a saber cual sea la tupla que deja el puesto $n + 1$ desocupado en su desenlace.
Por tanto, hay $(n+1)^n / (n+1) = (n+1)^{n-1}$ funciones de parqueo de longitud $n$.

Ahora, considere alguna $n$-tupla arbitraria $\mathbf{x} = (x_1, x_2, \ldots, x_n) \in [n]^n$.
?`Es posible determinar si es una función de parqueo de longitud $n$ sin explícitamente considerar su desenlace?
Una caracterización clásica basada en desigualdades permite esto precisamente.
Sea $\mathbf{x}' = (x_1', x_2', \ldots, x_n')$ el reordenamiento de $\mathbf{x}$ en orden creciente.
Entonces, lo siguiente se cumple:
\begin{align}
\label{eq: characterization}
    \parbox{0.375\linewidth}{
        \centering 
        $\mathbf{x}$ es una función de parqueo de longitud $n$
    } 
    &&
    \iff
    &
    \parbox{0.375\linewidth}{
        \centering 
        $x_i' \leq i$ para todo $1 \leq i \leq n$
    }
\end{align}
Esta caracterización se puede verificar de la siguiente manera.
Si existe alguna $i$ tal que $x_i' > i$, entonces al menos $n - i + 1$ carros intentarían parquear en los $n - i$  últimos puestos $s_{i + 1}, s_{i + 2}, \ldots, s_n$.
Estos serían más carros que puestos disponibles, por lo que al menos uno de estos carros no lograría parquear.
Por el contrario, si algún carro no logra parquear, entonces a fin de cuentas deberá existir un primer puesto disponible $s_i$, lo cual solamente puede ocurrir si $x_i' > i$.

Observe que \eqref{eq: characterization} tiene una implicación extraordinaria: las funciones de parqueo son invariantes bajo la acción del grupo de simetría $\mathfrak{S}_n$, el mismo que permuta sus índices.
En particular, el conjunto de funciones de parqueo se obtiene mediante los reordenamientos de las funciones de parqueo en orden creciente. 
Por ejemplo, $(1,3,1,1)$, $(1,1,3,1)$, $(1,1,1,3)$, y $(3,1,1,1)$ son todas funciones de parqueo de longitud $4$, cada cual con un desenlace distinto, mientras que ningún reordenamiento de $(1,3,3,4)$ podrá ser una función de parqueo de longitud $4$.

Las funciones de parqueo también están estrechamente relacionadas a los números de Catalan (OEIS \href{https://oeis.org/A000108}{A000108}), los cuales surgen de la relación de recurrencia
\begin{equation}
\label{eq: catalan}
    C_{n} 
    = \sum_{i=1}^{n} C_{i-1} \cdot C_{n - i}
\end{equation}
para $n \geq 1$ con $C_0 = 1$.
En particular, el conjunto de funciones de parqueo en orden creciente de longitud $n$ está enumerado por \eqref{eq: catalan}.
Para verificar esta relación, suponga que construye una función de parqueo en orden creciente de longitud $n$, denotada $\mathbf{x}' = (x_1', x_2', \ldots, x_n') \in [n]^n$, con su elección de $1 \leq i \leq n$ para el índice mas grande que satisfaga $x_i' = i$.
Para cualquiera sea su elección de $i$, la $(i-1)$-tupla $(x_1', x_2', \ldots, x_{i-1}')$ debe ser una función de parqueo de longitud $i-1$, la $(n-i)$-tupla $(x_{i+1}' - i + 1, x_{i+2}' - i + 1, \ldots, x_n' - i + 1)$ debe ser una función de parqueo de longitud $n-i$, y por inducción hay $C_{i-1} \cdot C_{n - i}$ posibilidades.

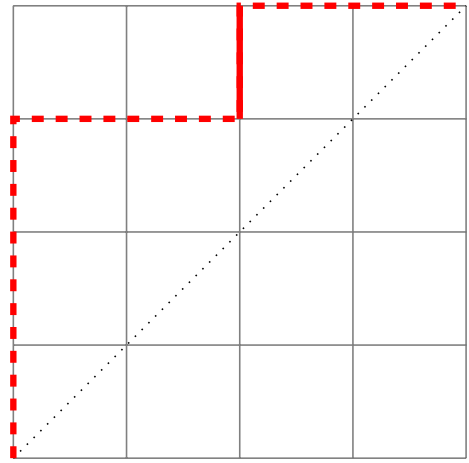
\begin{figure}[ht]
\centering
\resizebox{0.8\linewidth}{!}{
\begin{tikzpicture}
    \draw[step=1cm,gray,thin] (0,0) grid (4,4);
    \draw[thin, dotted] (0,0)--(4,4);
    \draw[ultra thick,dashed,red] (0,0)--(0,1)--(0,2)--(0,3)--(1,3)--(2,3)--(2,4)--(3,4)--(4,4);
    \draw[ultra thick,red] (2,3)--(2,4);
\end{tikzpicture}
}
    \caption{
        El camino de Dyck $NNNEENEE$ de longitud $8$, el cual está resaltado por la línea punteada roja en negrita, corresponde a la función de parqueo en orden creciente $\mathbf{x}' = (x_1', x_2', x_3', x_4') = (1,1,1,3)$ de longitud $4$.
        Por ejemplo, el cuarto paso $N$ en el camino es el sexto paso en total, y está resaltado por la línea continua roja en negrita.
        Hay un total de dos pasos $E$ anteriores, por lo que $x_4' = 1 + 2 = 3$.
    }
    \label{fig: dyck path}
\end{figure}

En consecuencia, las funciones de parqueo en orden creciente están en biyección con la gran variedad de objetos de Catalan; refiérase a Stanley~\cite{stanley2015catalan} para un vistazo general.
Por ejemplo, las funciones de parqueo en orden creciente de longitud $n$ están en biyección con los caminos de Dyck de longitud $2n$; estos son caminos reticulares desde $(0,0)$ hasta $(n,n)$ utilizando únicamente pasos hacia el norte $(1,0)$ y hacia el este $(0,1)$, denotados ``$N$'' y ``$E$'' respectivamente, y que no cruzan por debajo de la diagonal principal.
Armstrong, Loehr, y Warrington~\cite{armstrong2016rational}*{Sección~2.2} describen la siguiente biyección: dado un camino de Dyck de longitud $2n$, obtenga una función de parqueo en orden creciente $\mathbf{x}' = (x_1', x_2', \ldots, x_n')$ de longitud $n$ tras asignar a $x_i'$ uno más el número total de pasos $E$ que aparecen antes del $i$-ésimo paso $N$, para toda $1 \leq i \leq n$.
La Figura~\ref{fig: dyck path} demuestra esta construcción utilizando la función de parqueo en orden creciente $(1,1,1,3)$; el reordenamiento creciente del ejemplo de la Figura~\ref{fig: parking example}.

Note que si $\mathbf{x}' = (x_1', x_2', \ldots, x_n') \in [n]^n$ es una función de parqueo en orden creciente, entonces el $i$-ésimo carro $c_i$ (con preferencia $x_i'$) parquea en el $i$-ésimo puesto $s_i$.
Por lo tanto, el desplazamiento total de una función de parqueo en orden creciente $\mathbf{x}'$ es dado por
\begin{equation}
\label{eq: displacement}
    \sum_{i=1}^n i - x_i' = \frac{n(n+1)}{2} - \sum_{i=1}^n x_i'.
\end{equation}
Esta estadística tiene otras interpretaciones combinatoriales.
Por ejemplo, el número de celdas enteras entre un camino de Dyck y la diagonal principal, que en el caso de la Figura~\ref{fig: dyck path} es $4$, es igual al desplazamiento total de una función de parqueo en orden creciente.
Por último, note que la dependencia de \eqref{eq: displacement} en $\mathbf{x}'$ se reduce a la suma de sus términos.
!`Esto implica que el desplazamiento total de una función de parqueo también se invariante bajo reordenamientos!

\begin{figure}[ht]
    \centering
    \includegraphics[width=\linewidth]{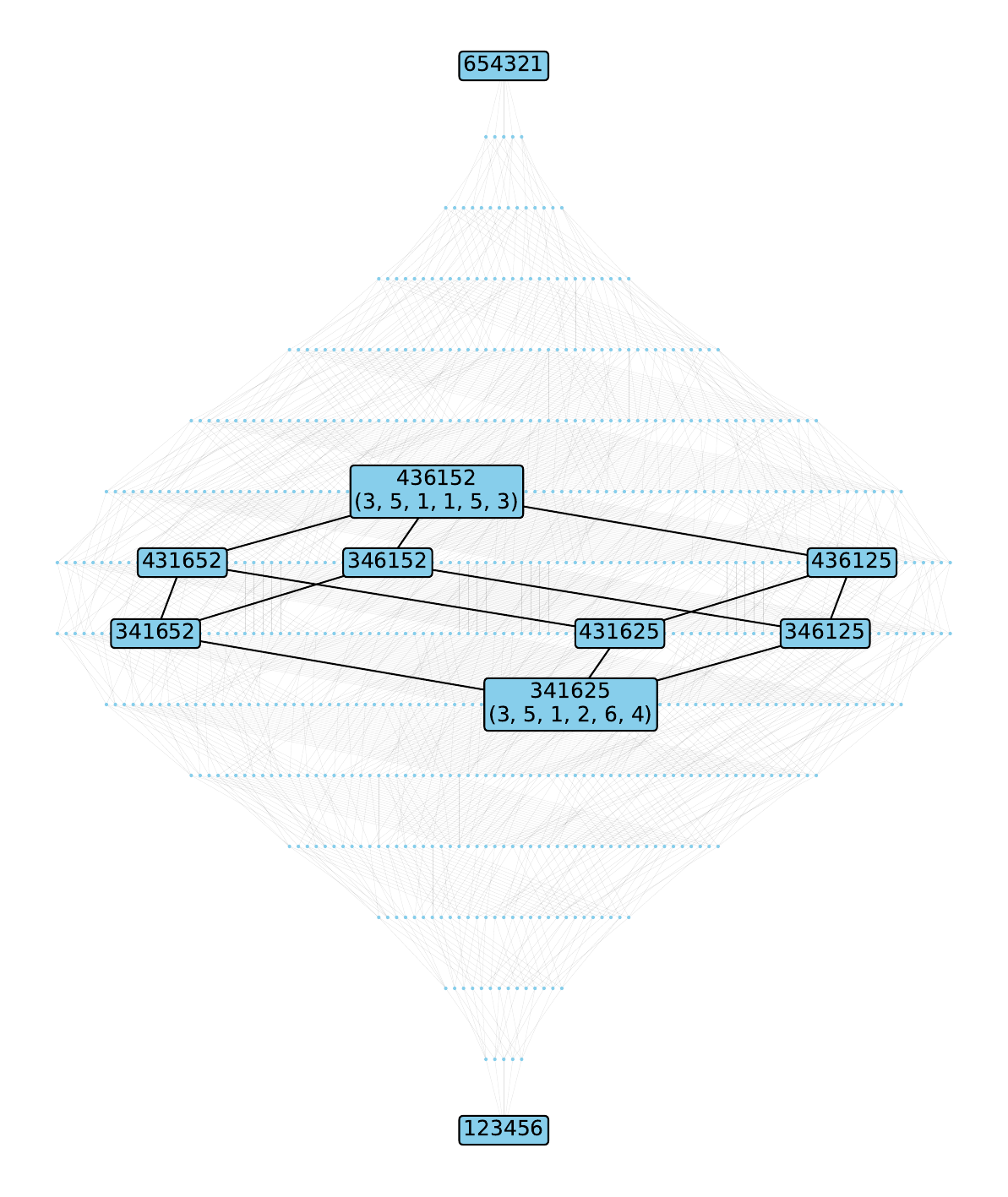}   
    \caption{
        Conjunto parcialmente ordenado de $\mathfrak{S}_6$ con un intervalo de Boole de rango $3$ resaltado, con elemento mínimo $341625$ y elemento máximo $436152$ escritos en notación de una línea.
        El desenlace de $(3,5,1,2,6,4)$, una función de parqueo de longitud $6$, es $(3, 4, 1, 6, 2, 5)$ y corresponde al elemento mínimo del intervalo.
        El desenlace de $(3,5,1,1,5,3)$, otra función de parqueo de longitud $6$, es así mismo $(3, 4, 1, 6, 2, 5)$ y corresponde al elemento mínimo del intervalo.
        De hecho, $(3,5,1,1,5,3)$ es una ``clasificación unitaria de Fubini con $3$ rangos distintos'' y corresponde a un intervalo de Boole de rango $6-3 = 3$ (i.e., un cubo), mientras que $(3,5,1,2,6,4)$ es una ``clasificación unitaria de Fubini con $6$ rangos distintos'' y corresponde a un intervalo de Boole de rango $6 - 6 =0$ (i.e., un vértice).
        Refiérase a Elder et al.~\cite{elder2024boolean} para más detalles.
    }
    \label{fig: weak}
\end{figure}

A pesar de que aparecieron en la literatura por primera vez hace más de seis décadas, la combinatoria de las funciones de parqueo se mantiene como una emocionante rama de investigación.
Trabajos recientes van desde el estudio de estadísticas discretas tal como el desplazamiento~\cites{kenyon2023parking,elder2023cost} o el número de carros suertudos~\cites{gessel2006refinement,slavik2023lucky,stanley2023some}, la introducción de nuevas variantes y/o generalizaciones de las reglas del ``experimento de parqueo'' clásico (refiérase a Carlson et al.~\cite{carlson2021parking} para un vistazo del sin fin de posibilidades al estilo de  ``elija su propia aventura''), aspectos poliedrales~\cites{amanbayeva2022convex,hanada2023generalized}, y sorprendentes conexiones con objetos que aparentan no estar relacionados~\cites{aguillon2023parking,harris2024lucky}, por mencionar solamente unos pocos.
En un trabajo reciente, mis colaboradores y yo utilizamos un subconjunto de las funciones de parqueo que llamamos \emph{clasificaciones unitarias de Fubini}, y en particular su desenlace, para caracterizar y contar los intervalos de Boole de rango $k$ en un conjunto parcialmente ordenado de $\mathfrak{S}_n$~\cite{elder2024boolean}.
La Figura~\ref{fig: weak} demuestra nuestra construcción.

\bibliography{bib}